\documentclass[12pt]{article}

\usepackage{euscript,amsmath, amssymb, amsfonts}
\usepackage{color}

\newcommand{\HH}{\mathcal H }
\newcommand{\HG}{$\mathcal H$}

\newcommand{\G}{{\mathbb C}[I_2(2m+1)]}
\newcommand{\spur}{sp^{\varkappa} }

\newcommand{\F}{F_0^{sp^{\varkappa}} }

\newcommand{\m}{\mathfrak {s} }

\newcommand{\Lo}{\mu L_0 }

\newcommand{\be}{\begin{equation}}
\newcommand{\ee}{\end{equation}}
\newcommand{\bee}{\begin{eqnarray}}
\newcommand{\eee}{\end{eqnarray}}
\newcommand\nn{\nonumber \\}

\newcommand\defeq{:=}

\newcounter{theorem}
\newcommand{\theorem}{\par\refstepcounter{theorem}
           {\bf Theorem \arabic{section}.%
           \arabic{theorem}. }}
\renewcommand\thetheorem{\thesection.%
   \arabic{theorem}}
\makeatletter \@addtoreset{theorem}{section}

\newcounter{corollary}
\newcommand{\corollary}{\par\refstepcounter{theorem}
           {\bf Corollary \arabic{section}.%
           \arabic{theorem}. }}

\makeatletter \@addtoreset{corollary}{section}

\newcounter{lemma}

\makeatletter \@addtoreset{lemma}{section}

\newcounter{proposition}
\newcommand{\proposition}{\par\refstepcounter{theorem}
           {\bf Proposition \arabic{section}.%
           \arabic{theorem}. }}
\renewcommand\theproposition{\thesection.%
  \arabic{theorem}}
\makeatletter \@addtoreset{proposition}{section}

\newcounter{conjecture}
\newcommand{\conjecture}{\par\refstepcounter{theorem}
           {\bf Conjecture \arabic{section}.%
           \arabic{theorem}. }}

\makeatletter \@addtoreset{conjecture}{section}

\newcounter{remark}
\newcommand{\remark}{\par\refstepcounter{theorem}
           {\bf Remark \arabic{section}.%
           \arabic{theorem}. }}
\renewcommand\theremark{\thesection.%
     \arabic{theorem}}
\makeatletter \@addtoreset{remark}{section}

\newcounter{definition}
\newcommand{\definition}{\par\refstepcounter{theorem}
           {\bf Definition \arabic{section}.%
           \arabic{theorem}. }}
\renewcommand\thedefinition{\thesection.%
    \arabic{theorem}}
\makeatletter \@addtoreset{definition}{section}

\newenvironment{proof}[1][Proof]{\noindent\textsf{#1.\ }}
{\hfill {\small $\square$}}

\makeatletter \@addtoreset{equation}{section}
\def\theequation{\thesection.\arabic{equation}}

\begin{document}


\sloppy \title
 {
Ideals generated by traces or by supertraces in the symplectic
reflection algebra $H_{1,\nu}(I_2(2m+1))$ II
}

\author
 {
  I.A. Batalin%
\thanks{ I.E. Tamm Department of Theoretical Physics,
          P.N. Lebedev Physical Institute, RAS
          119991, Leninsky prosp., 53, Moscow, Russia.}
\thanks{E-mail: batalin@lpi.ru}
 ,
 S.E. Konstein$^*$%
\thanks{E-mail: konstein@lpi.ru  (corresponding author) }
 ,
 I.V. Tyutin$^*$%
\thanks{E-mail: tyutin@lpi.ru}
}

\date{
}

\maketitle
\thispagestyle{empty}

\begin{abstract}

The algebra $\mathcal H:= H_{1,\nu}(I_2(2m+1))$ of observables of the Calogero model based on the root system
$I_2(2m+1)$ has an $m$-dimensional space of traces and an $(m+1)$-dimensional space of supertraces.
In the preceding paper we found
all values of the parameter $\nu$
 for which either the space of traces contains a~degenerate nonzero trace $tr_{\nu}$ or
the space of supertraces contains a~degenerate nonzero supertrace $str_{\nu}$
and, as a~consequence, the algebra $\mathcal H$ has
two-sided ideals: one consisting of all vectors
in the kernel of the form
$B_{tr_{\nu}}(x,y)=tr_{\nu}(xy)$ or another consisting of all vectors in the kernel of
the form
$B_{str_{\nu}}(x,y)=str_{\nu}(xy)$.
We noticed that if $\nu=\frac z {2m+1}$, where $z\in \mathbb Z \setminus (2m+1) \mathbb Z$, then there exist both a
degenerate trace and a~degenerate supertrace on $\mathcal H$.

Here we prove that the ideals determined by these degenerate
forms coincide.

\end{abstract}




\section{Introduction}
This paper is a continuation of \cite{KT7}; we advise the reader to recall \cite{KT7}.

\subsection{Definitions}

Let $\mathcal A$ be an associative $\mathbb Z_2$-graded algebra with unit;
let  $\varepsilon$ denote its parity. All expressions of linear algebra are given
for homogenous elements only and are supposed to be extended to inhomogeneous elements via
linearity.

A linear complex-valued function $tr$ on $\mathcal A$ is called a~\emph{trace}
if $tr(fg-gf)=0$ for all $f,g\in \mathcal A$.
A linear complex-valued function $str$ on $\mathcal A$ is called a~\emph{supertrace}
if ${str(fg-(-1)^{\varepsilon(f)\varepsilon(g)}gf)=0}$ for all $f,g\in \mathcal A$.
These two definitions can be unified as follows.

Let ${\varkappa=\pm 1}$.
A linear complex-valued function  $\spur$ on $\mathcal A$ is called \emph{$\varkappa$-trace}
if ${\spur(fg-\varkappa^{\varepsilon(f)\varepsilon(g)}gf)=0}$ for all $f,g\in \mathcal A$.

Each nonzero  $\varkappa$-trace $\spur$ defines the nonzero
symmetric\footnote{Initially,
we used the term ``(super)symmetric bilinear form"
currently used by many, e.g., in the paper \cite{BKLS},
even in its title. However, in a recent preprint \cite{KLLS},
it is explained that the supersymmetry
$B(v,w)=(-1)^{p(v)p(w)}B(w,v)$ is related with
the isomorphism $V\otimes W\simeq W\otimes V$
of superspaces and has nothing to do with the (anti)symmetry of the bilinear form $B$ on $V=W$.}
bilinear form $B_{\spur}(f,g):=\spur(fg)$.

If $B_{\spur}$ is degenerate, then the set of the vectors of
its kernel is a~proper ideal in $\mathcal A$. We say that the
 $\varkappa$-trace $\spur$ is \emph{degenerate} if the bilinear form
$B_{\spur}$ is degenerate.

\subsection{The goal and structure of the paper}

The simplicity (or, alternatively, existence of ideals) of
Symplectic Reflection Algebras or, briefly, SRA (for definition, see \cite{sra})
was investigated in a~number of papers, see, e.g.,  \cite{BG}, \cite{IL}.
In particular, it is shown that all SRA $H_{1,\nu}(G)$ with $\nu=0$
are simple (see \cite{pass}, \cite{BG}).

It follows from \cite{KT} and \cite{stek} that an associative algebra of observables of
the Calogero model with harmonic term in the potential and with coupling constant $\nu$
based on the root system $I_2(2m+1)$ (this algebra is SRA denoted $H_{1,\nu}(I_2(2m+1))$)
has an $m$-dimensional space of  traces and
an $(m+1)$-dimensional space of supertraces.

We say that the parameter $\nu$ is \emph{singular}, if the algebra $H_{1,\nu}(I_2(n))$ has
a~degenerate trace or a~degenerate supertrace.

In \cite{KT7}, we found
all singular values of $\nu$
for the algebras $\mathcal H := H_{1,\nu}(I_2(n))$ in the case of $n$ odd ($n=2m+1$) and found the
corresponding degenerate traces and supertraces;
the result is formulated in Theorem \ref{th2}.

We
noticed that if $\nu=\frac z {2m+1}$, where $z\in \mathbb Z \setminus (2m+1) \mathbb Z$,
then there exist both a
degenerate trace  and a~degenerate supertrace on $H$.

Denote this degenerate trace by $tr_z$ and the degenerate supertrace by $str_z$

Theorem \ref{th2} proved in \cite{KT7} implies that if $z\in \mathbb Z \setminus n\mathbb Z$, then

(i) the trace given by the formula (\ref{913})  in \cite{KT7} is degenerate and generates the ideal
$\mathcal I_{tr_z}$ consisting of all the vectors
in the kernel of the degenerate form $B_{tr_z}(x,y)=tr_z(xy)$, 

(ii) 
the supertrace (\ref{914}) is degenerate and generates the ideal $\mathcal I_{str_z}$ consisting of all the vectors in the kernel of
the degenerate form ${B_{str_z}(x,y)=str_z(xy)}$.


The goal of this paper is
Theorem \ref{conj9.1}, which proves
{\conjecture\label{con1} (\cite[Conjecture 9.1]{KT7})
$\mathcal I_{tr_z}=\mathcal I_{str_z}$.
}

\vskip 3mm

In Sections \ref{i2}--\ref{spec} we recall the necessary definitions
and
preliminary facts.

\section{The group $I_2(2m+1)$}\label{i2}

Hereafter in this paper, $n=2m+1$.

\definition
The group $I_2(n)$ is a~finite subgroup of the orthogonal group $O(2,\mathbb R)$ generated
by the root system $I_2(n)$.

The group $I_2(n)$ is the symmetry group of a~flat regular $n$-gon;
$I_2(n)$
consists of $n$
reflections $R_k$ and $n$ rotations $S_k$, where $k=0,1,...,2m$. We consider the indices $k$ as
integers modulo $n$.

These elements ($R_k$ and  $S_k$ for all $k$)
satisfy the relations
\be
\nonumber
R_k R_l = S_{k-l},\qquad
S_k S_l = S_{k+l},\qquad
R_k S_l = R_{k-l},\qquad
S_k R_l = R_{k+l}.
\ee
The element $S_0$
is the unit in the group $I_2(n)$.
Obviously, since $n$ is odd,
all the reflections $R_k$
are in the same conjugacy class.

The rotations $S_k$ and $S_l$ constitute a~
conjugacy class if
$k+l=n$.

{
Let
\be
           \nonumber
\lambda:= \exp\left(\frac {2\pi i}{n}\right).
\ee

Let
\be\label{defg}
G:={\mathbb C}[I_2(n)]
\ee
be the group algebra of the group $I_2(n)$.
In $G$, it is convenient to introduce the following basis
\bee
\nonumber
&& L_p:=\frac 1 n \sum_{k=0}^{n-1} \lambda^{kp} R_k,\qquad
Q_p:=\frac 1 n \sum_{k=0}^{n-1} \lambda^{-kp} S_k.
\eee


\section{Symplectic reflection algebra
$H_{1,\nu}(I_2(2m+1))$}\label{srai2}


{\definition
\emph{The symplectic reflection algebra
$\mathcal H :=H_{1,\nu}(I_2(2m+1))$} is the associative algebra
of  polynomials in the noncommuting elements $a^{\alpha}$ and $b^{\alpha}$, where $\alpha=0,1$,
with coefficients in $G$ (see Eq. (\ref{defg})), satisfying the relations%
}
\bee
\nonumber
&L_p a^\alpha = -   b^\alpha L_{p+1} ,\qquad
&L_p b^\alpha = -   a^\alpha L_{p-1} ,\nn
\label{QH}
&Q_p a^\alpha = a^\alpha Q_{p+1} ,\qquad
&Q_p b^\alpha = b^\alpha Q_{p-1} ,\\
&L_k L_l =  \delta_{k+l} Q_l, \qquad &L_k Q_l =  \delta_{k-l} L_l,\nn
&Q_k L_l =  \delta_{k+l} L_l, \qquad & Q_k Q_l =  \delta_{k-l} Q_l,
\nonumber
\eee
\bee
\left[a^\alpha,\,b^\beta\right]&=& \varepsilon^{\alpha\beta} \hfill
\left( 1+ \mu  L_0 \right), \nn
\left[a^\alpha,\,a^\beta\right]&=&\varepsilon^{\alpha\beta}
\mu   L_1 ,\nn
\left[b^\alpha,\,b^\beta\right]&=&\varepsilon^{\alpha\beta}
\mu   L_{-1},
     \nonumber
     \eee
where $\delta_k\defeq \delta_{k0}$, and
$\varepsilon^{\alpha\beta}$ is the skew-symmetric tensor normalized so that $\varepsilon^{0\,1}=1$, and
 \be
 \nonumber
 \text{~~}\mu:=n {\nu}.
\ee

Defining the parity in \HG  \ by setting
\be     \nonumber
\varepsilon(a^\alpha)=\varepsilon(b^\alpha)=1,\ \ \ \varepsilon(R_k)=\varepsilon(S_k)=0,
\ee
we turn this algebra into a~superalgebra.

The algebra $H_{1,\nu}(I_2(2m+1))$ depends on one complex parameter $\nu$.


\section{Subalgebra of singlets}\label{sl2}

Consider the elements%
\footnote
{Here the brackets $\{\cdot ,\cdot\}$ denote anticommutator.
}
$
T^{\alpha\beta}:=\frac 1 2 (\{a^\alpha,\, b^\beta \}+\{b^\alpha,\, a^\beta \})
$
of the algebra \HG,
\ and the inner derivations of  \HG they generate:
\be     \nonumber
D^{\alpha\beta}:\ \ f\mapsto \left[f,T^{\alpha\beta}\right] \ \
\text{for any } f\in\HH.
\ee
It is easy to verify that the linear span of these derivations
is a~ Lie algebra isomorphic to $sl_2$.

\definition
A \emph{singlet} is any element $f\in \HH$ such that $[f,\, T^{\alpha\beta}]=0$ for all $\alpha,\beta$.
The subalgebra
$H^0 \subset \HH$ consisting of all singlets of the algebra
\HG \ is called  the \emph{subalgebra of singlets}.

One can consider the algebra \HG \ as an $sl_2$-module and decompose it into
the direct sum of irreducible submodules.

Observe, that any $\varkappa$-trace is identically zero on all irreducible
$sl_2$-submodules of \HG,  \ except for singlets.

Let the skew-symmetric tensor $\varepsilon_{\alpha\beta}$ be normalized so that $\varepsilon_{0\,1}=1$
and so $\sum_{\alpha} \varepsilon_{\alpha\beta} \varepsilon^{\alpha\gamma}=\delta_{\beta}^{\gamma}$.
We set
\bee     \nonumber
\m:=\frac 1 {4i} \sum_{\alpha,\beta=0,1}
\left(\{a^\alpha,\, b^\beta \}-\{b^\alpha,\, a^\beta \}\right)\varepsilon_{\alpha\beta}.
\eee

{\proposition \label{1.3.7} (\cite[Proposition 4.2]{KT7})
{\it The subalgebra of singlets $H^0$
is the algebra of polynomials in
the element $\m$ with coefficients in the group algebra  $\G$.}
}

The commutation relations of the singlet $\m$ with generators of
the algebra \HG  \ have the form:
\bee
[\m ,\,Q_p] = [\m ,\,S_k] = [T^{\alpha\beta},\,\m]=0,     \nn
\m L_p =-L_p\m, \qquad \m R_k =-R_k\m,          \nonumber  \nn
(\m -i\mu L_0 )a^\alpha = a^\alpha (\m + i+i \mu L_0).       \nonumber
\eee

\section{The form of ideals in \HG\  and in $H^0$ }\label{id}

{\theorem\label{i0} (\cite[Theorem 4.3]{KT7})
Let $\mathcal I$ be a~proper ideal in the algebra \HG, and
$\mathcal I_0 :=\mathcal I \bigcap H^0$.
Then, there exist  nonzero polynomials
$\phi_k^0\in \mathbb C[\m]$, where $k=0,...,n-1$, such that
$\mathcal I_0$ is the span over $\mathbb C[\m]$ of the elements
\bee               \nonumber
\phi_k^0 (\m)Q_k, \qquad \ \phi_{n-k}^0 L_k, \text{ where } k=0,\,...\,,\,n-1
\ \text{and }\ \phi_n^0:=\phi_0^0.
\eee
}

{\proposition\label{ne0} (\cite[Proposition 4.4]{KT7})
If $\mathcal I\subset \HH$ is a~proper ideal,
then
$\mathcal I_0 =\mathcal I \bigcap H^0$ is a~proper ideal in
$H^0$.
}

\definition
For each $p=0, ... , 2m$, we define the ideals $\mathcal J_p$ and $\mathcal J^p$
in the algebra $\mathbb C[\m]$
by setting
\be     \nonumber
\mathcal J_p := \{f\in \mathbb C[\m] \mid \ f(\m)Q_p \in \mathcal I    \},
\qquad \mathcal J^p := \{f\in \mathbb C[\m] \mid \ f(\m)L_p \in \mathcal I    \}.
\ee

{\proposition\label{8.13} (\cite[Proposition 4.7]{KT7}).
 We have
$\mathcal J_p=\mathcal J^{-p}$.
}

{\proposition\label{9.4} (\cite[Proposition 4.8]{KT7}).
 We have $\mathcal J_p\ne 0$ for any $p=0,...,2m$.
}

Since $\mathbb C[\m]$ is a~principal ideal ring, we have the following statement:

 {\corollary \label{col}   For any
$p=0,...,2m$, there exists a~nonzero polynomial $\phi_p^0\in \mathbb C[\m]$
such that $\mathcal J_p =\phi_p^0 \mathbb C[\m]$.
}

Theorem \ref{i0} evidently follows from Corollary \ref{col}.


\section{Generating functions of $\varkappa$-traces}\label{genfuni2}

For each $\varkappa$-trace $\spur$ on \HG, one can define the following
set of generating functions  which allow one to calculate the $\varkappa$-trace
of arbitrary element in $H^0$ via finding the values of the derivatives of these functions with respect to parameter $t$ at zero:
\bee\label{1.33}
&& F_p^{\spur}(t):=\spur(\exp(t(\m-i\Lo))Q_p),
  \\
&& \Psi_p^{\spur}(t):=\spur(\exp(t\m)L_p),   \text{~~where $p=0,...,2m$.}  \nonumber
\eee

Since $L_0 Q_p=0$ for any $p\ne 0$,
it follows from the definition
~(\ref{1.33}) that
\bee\nonumber
F_p^{\spur}(t)&=& \spur(\exp(t \m)Q_p) \text{ if } p\ne 0,
\nn
F_0^{\spur}(t)&=&\spur(\exp(t(\m-i\mu L_0))Q_0).
\nonumber
\eee

It is easy to find $\Psi_p^{\spur}$ for $p\ne 0$.
Since $\m L_q =-L_q\m$ for any $q=0,...,2m$,
we have
\be\label{psiq}
\Psi_q^{\spur}(t)=\spur(\exp(t\m)L_q)=\spur(L_q).
\ee
Next, since $\spur(R_k)$ does not depend on $k$, we have
$\spur(L_p)=0$ for any $p\ne 0$ and
\be\label{psi0}
\Psi_p^{\spur}(t)\equiv 0 \text{ for any } p \ne 0.
\ee

The value of $\spur(L_0)$ will be calculated later, in Section \ref{GAi2-sp}.

We consider also the functions
\[
\Phi_p^{\spur}(t):=\spur(\exp(t(\m+i\Lo))Q_p).
\]
It is easily verified, by expanding the exponential in a~series, that
these functions are related with the functions $F_p^{\spur}$
by the formula
\be     \nonumber
\Phi_p^{\spur}(t)=F_p^{\spur}(t)+2i\Delta_p^{\spur}(t), \\
\text{where } \Delta_p(t)^{\spur}=\delta_p \sin(\mu t)\spur(L_0).
     \nonumber
\ee

The form of generating functions is related with (non)degeneracy of the form
$B_{\spur}$  as described in Proposition \ref{1.7.8} below.


\section{Degeneracy conditions for the $\varkappa$-trace}\label{degcond}

{\proposition\label{1.7.8} (\cite[Proposition 6.1]{KT7}).
The $\varkappa$-trace on the algebra \HG \ is degenerate  if and only if the
generating functions $F_p^{\spur}$ defined by  formula \eqref{1.33}
have the following form
\bee\label{wid1}
F_p^{\spur}(t)=\sum_{j=1}^{j_p} \exp(t \omega_{j,p})\varphi_{j,p}(t),
\eee
where $\omega_{j,p}\in \mathbb C$ and $\varphi_{j,p}\in \mathbb C [t]$ might depend on $\varkappa$.
}


\section{Equations for the generating functions $F_p^{\spur}$}\label{genfun}

In \cite[Eq. (7.1)]{KT7}, the following system of differential equations for the generating functions is obtained:
\bee\label{eqgenfun}
\frac d {dt} F_p^{\spur} - \varkappa e^{it}\frac d {dt}F_{p+1}^{\spur}=iF_p^{\spur} +
\varkappa ie^{it}F_{p+1}^{\spur}+2 \varkappa i\frac d {dt}\left(e^{it}\Delta_{p+1}^{\spur} \right).
\eee

The initial conditions for this system are:
$$
F_p^{\spur}(0) = \spur(Q_p).
$$

To solve the system (\ref{eqgenfun}), we consider its Fourier transform.
Let
\bee     \nonumber
&& \lambda:=e^{2\pi i/(2m+1)}, \\
&& G_k^{\spur}:= \sum_{p=0}^{2m} \lambda^{kp}F_p^{\spur}, \text{ where } k=0,...,2m,
\label{fur}
\\
&& \widetilde{\Delta}_k^{\spur}:=\sum_{p=0}^{2m} \lambda^{kp}\Delta_{p+1}^{\spur}=
\lambda^{-k}\left( \sin(\mu t)\spur(L_0)    \right),
 \text{ where } k=0,...,2m.
     \nonumber
\eee
For the functions $G_k^{\spur}$, we then obtain the equations
\be\label{Geq}
\frac d {dt} G_k^{\spur} = i\frac {\lambda^k+ \varkappa e^{it}}{\lambda^k- \varkappa e^{it}} G_k^{\spur} +
\frac {2i \varkappa \lambda^k }{\lambda^k- \varkappa e^{it}} \frac d {dt}\left(e^{it}\widetilde{\Delta}_k^{\spur}\right)
\ee
with the initial conditions
\be     \label{init}
G_k^{\spur}(0)=\spur(S_k).
\ee

We choose the following form of the solution of the system (\ref{Geq}):
\bee\label{GK}
G_k^{\spur}(t)= \frac {\varkappa e^{it}} {(\varkappa e^{it}-\lambda^k)^2} \lambda^k g_k^{\spur}(t),
\eee
where
\bee\label{1.47g}
\begin{array}{ll}
g_k^{\spur}(t)&=\left(\frac{2}{\mu}(\cos (t\mu) -1)
+{2i}\lambda ^{-k}(\lambda ^{k}-\varkappa e^{it})\sin (t\mu)
\right)\spur(L_0)\\
&+\varkappa\lambda ^{-k}(\varkappa -\lambda ^{k})^{2} \spur(S_k).
\end{array}
\eee

Evidently, this solution satisfies the initial condition (\ref{init}) for each $\varkappa$
and $k$, except for the case where $\varkappa=+1$ and $k=0$.


If $\varkappa = +1$ and $k=0$, then the expression
~(\ref{GK})
for $G_0^{tr}$ has a~removable singularity at $t=0$.
In this case, instead of the condition ${G_0^{tr}(0)=tr(S_0)}$
we consider the condition $\lim_{t \rightarrow 0} G_0^{tr}(t) = tr(S_0)$.

When $\varkappa=+1$ the solution (\ref{GK}) -- (\ref{1.47g}) gives
\be                 \nonumber
G_0^{tr}(t)= \frac { e^{it}} {( e^{it}-1)^2}
\left(\frac{2}{\mu}(\cos (t\mu) -1)
+{2i}(1- e^{it})\sin (t\mu)
\right)tr(L_0),
\ee
and one can easily see that
\be   \nonumber
\lim_{t \rightarrow 0} G_0^{tr}(t) = -\mu tr(L_0).
\ee

It is shown in Subsection \ref{GAi2-tr}
that if $\varkappa=+1$,
then
\be    \nonumber
tr(S_0) = -\mu tr(L_0)
\ee
for any trace $tr$ on \HG.

So, $G_0^{tr}(t)$ satisfies the initial conditions (\ref{init}) also.


In the case where $ \varkappa = -1$, the $\varkappa$-trace is a~supertrace (see \cite{KT}).
In this case, the $m+1$ values $str(S_k)=str(S_{2m+1-k})$ for $k=0,...,m$ completely define
the supertrace on \HG \ (see \cite{stek}).

In the case where $ \varkappa = +1$, the $\varkappa$-trace is a~trace (see \cite{KT}).
In this case, the $m$ values ${tr(S_k)=tr(S_{2m+1-k})}$ for $k=1,...,m$ completely define
the trace on \HG \ (see \cite{stek}).
The value $tr(S_0)$ linearly depends on parameters $tr(S_k)$, where $k=1,...,m$,
and this value is found in Subsection \ref{GAi2-tr} (see Eqs. (\ref{s0}) -- (\ref{9.6})).

\section{Values of the $\varkappa$-trace on $\G$}\label{GAi2-sp}

From \cite{KT7} we have
\be\label{rk}
\spur (R_k) = - \frac {2\mu} {2m+1}
\left(\frac
{1+\varkappa} 2 X^{tr}
+
\frac {1-\varkappa} 2 Y^{str}
\right),
\ee
where
\bee\label{92}
&& X^{tr}:= \sum_{r=1}^{2m}
\sin^2\left(\frac {\pi r}{2m+1}\right) tr (S_r),\\
&& Y^{str}:= \sum_{r=0}^{2m}                            \label{93}
\cos^2\left(\frac {\pi r}{2m+1}\right) str (S_r).
\eee

Below we consider these values for the traces and supertraces separately.

\subsection{Values of the traces ($\varkappa=+1$) on $\G$}\label{GAi2-tr}

The group $I_2(2m+1)$ has $m$ conjugacy
classes without the eigenvalue +1 in the spectrum:
$
\{S_{p},S_{2m+1-p}\}, \text { where } p=1,...,m.
$

By Theorem 2.3 in \cite{KT}, the values of the trace on these conjugacy classes
\be
\nonumber
s_k := tr(S_k),\ \ \text{ where } s_{2m+1-k}=s_k, \ \ k=1,...,m,
\ee
are arbitrary and
completely define the trace on the algebra \HG.
Therefore, the dimension of the space of traces is equal to $m$.

Further,
the group $I_2(2m+1)$ has one conjugacy class with one eigenvalue +1 in its spectrum:
$
\{R_{1},\,...\,,\, R_{2m+1}\}.
$
The value of $tr(R_k)$ is expressed via $s_k$ by formula (\ref{rk}).

Besides, the group $I_2(2m+1)$ has
one conjugacy class with two eigenvalues +1 in its spectrum:~$\{S_0\}$.

The traces on conjugacy classes with two eigenvalues +1 in the spectrum
is calculated in \cite{KT7} using Ground Level Conditions (for their definition, see \cite{KT}):
\bee
\label{s0}
tr(S_{0}) = 2 \nu^{2}(2m+1) X^{tr}.
\eee

We also note that
\be\label{9.6}
tr(L_0)=-\frac {2\mu} {2m+1} X^{tr},\quad
tr (L_p)=0 \text{ for } p\ne 0, \quad
tr(S_0)=-\mu tr(L_0).
\ee


\subsection{Values of the supertraces ($\varkappa=-1$) on $\G$}\label{GAi2-str}

The group $I_2(2m+1)$ has $m+1$ conjugacy
classes without the eigenvalue $-1$ in the spectrum:
\be     \nonumber
\{S_{0}\},\
\{S_{p},S_{2m+1-p}\}, \text { where } p=1,...,m.
\ee
By \cite[Theorem 2.3]{KT}, the values of the supertrace on these conjugacy classes
\be
\nonumber
u_k := str(S_k)=str(S_{2m+1-k}),\ \ \text{ where }  \ \ k=0,...,m,
\ee
are arbitrary parameters that
completely define the supertrace $str$ on the algebra \HG,
and therefore the dimension of the space of supertraces is equal to $m+1$.

Besides,
the group $I_2(2m+1)$ has one conjugacy class with one eigenvalue $-1$ in the spectrum:
$\{R_{1},\,...\,,\, R_{2m+1}\}$.

The supertraces of the conjugacy class with eigenvalue $-1$ in  its spectrum
are given by Eq. (\ref{rk}):
$str(R_{k}) = -2\nu Y^{str}$, where
$k=0,1,...,2m$, and
where $Y^{str}$ is defined by Eq (\ref{93}).


\section{Singular values of the parameter $\mu$}\label{spec}

 The solution Eq (\ref{GK})-(\ref{1.47g}) determines the generation functions of traces and supertraces
on $H^0$ for any trace and any supertrace on \HG\ . Generally speaking, $G_k^{\spur}$ is a~meromorphic
function on $t$, but if $\mu$ and $\spur$ are such that the form $B_{\spur}$ is degenerate,
then $G_k^{\spur}$ is an integer function on $t$ for each $k$. The complete list of such pairs of $\mu$
and $\spur$ is given in Theorem \ref{th2}. For these values of $\mu$ and $\spur$, the functions
$G_k^{\spur}$ are Laurent polynomials in $\exp(it)$.

{\theorem\label{th2} (\cite[Theorem 9.1]{KT7}).
Let $ m\in \mathbb Z$, where $m\geqslant 1$,
and $n=2m+1$.
Then

\textup{1)} The associative algebra $H_{1,\nu}(I_2(n))$
has a~$1$-parameter set of nonzero  traces $tr_z$ such that the symmetric invariant bilinear form
$B_{tr_z}(x,y)=tr_z(xy)$  is degenerate if and only if
$\nu= \frac z n$, where $z\in \mathbb Z \setminus n\mathbb Z$.
These traces are completely defined by their values at
$S_k$ for
$k=1,\dots , m$:
\be\label{913}
tr_z(S_k)= \frac \tau {n \sin^{2}\frac {\pi k}n}\left
(1 - \cos\frac {2\pi k z}{n}\right), \text{ ~~where } \tau\in \mathbb C,\ \tau\ne 0.
\ee
Here $\tau$ is an arbitrary parameter specifying the trace in 1-dimensional space of
 traces.

\textup{2)}  The associative superalgebra $H_{1,\nu}(I_2(n))$
has a~$1$-parameter set of  nonzero supertraces $str_z$ such that the symmetric invariant bilinear form
$B_{str_z}(x,y)=str_z(xy)$  is degenerate if
$\nu= \frac z n$, where $z\in \mathbb Z \setminus n\mathbb Z$.
These supertraces are completely defined by their values at
$S_k$ for
$k=0,\dots , m$:
\bee\label{914}
str_z(S_k)= \frac \tau {n \cos^{2}\frac {\pi k}n}
\left(1 - (-1)^{z}\cos\frac {2\pi k z}{n}\right), \text{ where } \tau\in \mathbb C,\ \tau\ne 0.
\eee
Here $\tau$ is an arbitrary parameter specifying the supertrace in 1-dimensional space of supertraces.

\textup{3)}  The associative superalgebra $H_{1,\nu}(I_2(n))$
has a~$1$-parameter set of nonzero supertraces $str_{1/2}$ such that the symmetric invariant bilinear form
$B_{str_{1/2}}(x,y)=str_{1/2}(xy)$  is degenerate if
$\nu= z + \frac 1 2$, where $z\in \mathbb Z$.
These supertraces are completely defined by their values at
$S_k$ for
$k=0,\dots , m$:
\be             \nonumber
str_{1/2}(S_k)= \frac \tau {n \cos^{2}\frac {\pi k}n}, \text{ where } \tau\in \mathbb C, \ \tau\ne 0.
\ee
Here $\tau$ is an arbitrary parameter specifying the supertrace in 1-dimensional space of supertraces.

\textup{4)}  For all other values of $\nu$, all nonzero traces and supertraces are nondegenerate.
}

\section{Generating functions $F_p^{\spur}$   for the degenerate $\varkappa$-trace}

Let $\mu \in\mathbb Z \setminus n\mathbb Z$. Substitute the solutions
(\ref{913})   for the case $\varkappa=+1$  and
(\ref{914})   for the case $\varkappa=-1$
to Eqs. (\ref{GK})--Eq. (\ref{1.47g}). We obtain
the formula for both values of $\varkappa$
\be\label{g2}
g_k^{\spur}=-\frac{4\tau}{n}\left[
\cos(t\mu)+i\mu\lambda^{-k}(\lambda^k-\varkappa e^{it})\sin(t\mu)
-\varkappa^\mu\cos\frac{2\pi k\mu}{n}
\right].
\ee
Introducing the new variable $y$ instead of $t$
\be\label{y}
y:=\varkappa e^{it}
\ee
we can rewrite Eq. (\ref{g2}) in the form
\be
                       \nonumber
g_k^{\spur}=-\frac{2\tau}{n}\varkappa^\mu
\left[
(y^\mu+y^{-\mu})+\mu\lambda^{-k}(\lambda^k-y)(y^\mu-y^{-\mu})
-2\cos\frac{2\pi k\mu}{n}\right]
\ee
and Eq. (\ref{GK}) in the form
\be\label{GK2}
G_k^{\spur}=\frac{\lambda^k y}{(y-\lambda^k)^2}g_k^{\spur} .
\ee
Now we see that $G_k^{\spur}$ are the Laurent polynomials  in $y$ with the highest degree $\leq |\mu|$
and the lowest degree $\geq 1-|\mu|$.

Note, that the expressions  (\ref{GK2}) are even functions of the parameter $\mu$,
so we can assume that $\mu$ is a~positive integer.

Let $\mu > 0$ in what follows.

Thus, $G_k^{\spur}$ can be expressed in the form
\be\label{GK3}
G_k^{\spur}=\varkappa^\mu \sum_{\ell=\mu}^{1-\mu}\beta_{\ell}^{k}y^{\ell},
\ee
where
the 
$\beta_{\ell}^{k}$
are constants not depending on $\varkappa$ and
not all of them
equal to zero.

Eq (\ref{GK3}) implies that
\be\label{p0}
\beta_{\mu}^k=\frac{2\tau\mu}n.
\ee

Further, Eq. (\ref{fur}) implies
\be
       \nonumber
F_p^{\spur}=\frac 1 n \sum_{k=0}^{2m}\lambda^{-kp}G_k^{\spur}
\ee
and the generating functions $F_p^{\spur}$ have the form
\be\label{Fform}
F_p^{\spur}=\varkappa^\mu \sum_{\ell=\mu}^{-\mu}\alpha_{\ell}^{p}y^{\ell},
\ee
 where the $\alpha_{\ell}^{k}$ are constants not depending on $\varkappa$.
Observe that $F_p^{\spur}$ can be equal to zero for some $p\ne 0$ (e.g., if $\mu=1$, then $F_p^{\spur}=0$ for
each $p\ne 0$),
but $F_0^{\spur}\ne 0$
since Eq (\ref{p0}) implies $\alpha^0_{\mu}=\frac{2\tau\mu}n\ne 0$. Eq (\ref{p0}) implies also that $\alpha^0_{-\mu}=0$.


\section{The generating function $\mathcal F^{\spur}$=$\spur\left(\exp(t\m )Q_0 \right)$
  for the degenerate $\varkappa$-trace}

 Let $\mu \in\mathbb Z \setminus n\mathbb Z$ and $\varkappa$-trace be defined by
Eq (\ref{913}) in the case $\varkappa=+1$  and
by Eq (\ref{914}) in the case $\varkappa=-1$.

In this section we introduce the function
\be
           \nonumber
\mathcal F^{\spur}:=\spur\left(\exp(t \frak s )Q_0 \right)
\ee
and express it via $F_0^{\spur}$.
{\proposition\label{p121}
$\mathcal F^{\spur}$ is an even function of $t$:
\be\label{121}
\mathcal F^{\spur}=\spur(\cosh(t\m)Q_0 .
\ee
}

Indeed, $\mathcal F^{\spur}=\spur(\cosh(t\m)Q_0+\sinh(t\m)Q_0)$ and
$\spur(\sinh(t\m)Q_0)=0$ since
\bee\nonumber
\spur(\sinh(t\m)Q_0)=\spur((\sinh(t\m)L_0) L_0)=\spur(L_0(\sinh(t\m)L_0))=\nn=
\spur((L_0\sinh(t\m))L_0)=\spur((-(\sinh(t\m)L_0))L_0)
=\spur(-\sinh(t\m)Q_0)\nonumber
\eee

Now, decompose $\F$:
\bee
\F&=&\spur \left(e^{t(\m -i\mu L_0)} Q_0\right)=F_{even}+F_{odd}  \text{ where } \nn
F_{even}&=& \spur \left(\sum_{s=0}^{\infty}\frac 1 {(2s)!}(t(\m -i\mu L_0))^{2s} Q_0\right)
=\spur \left(\sum_{s=0}^{\infty}\frac 1 {(2s)!}t^{2s}(\m^2 -\mu^2)^{s} Q_0\right) ,                \label{even}
\\
F_{odd}&=& \spur \left(\sum_{s=0}^{\infty}\frac 1 {(2s+1)!}(t(\m -i\mu L_0))^{2s+1} Q_0\right)=  \nn
       &=&\spur \left(\sum_{s=0}^{\infty}\frac 1 {(2s+1)!}t^{2s+1}(\m^2 -\mu^2)^{s}
       (\m -i\mu L_0)Q_0\right)= \nn
       &=&
        \spur \left(\sum_{s=0}^{\infty}\frac 1 {(2s+1)!}t^{2s+1}(\m^2 -\mu^2)^{s}
       ( -i\mu L_0)Q_0\right) = \nn
&=&
        \sum_{s=0}^{\infty}\frac 1 {(2s+1)!}t^{2s+1}(-\mu^2)^{s}
       ( -i\mu)\spur L_0 =                       \nn
  &=&     \sinh (-i\mu)\spur L_0 =
       -\frac {\varkappa^\mu}{2}(y^{\mu}-y^{-\mu}) \spur L_0     .                 \label{odd}
\eee
Eq. (\ref{Fform}) implies that
\be\label{odd2}
F_{odd}= \frac {\varkappa^\mu} 2
\left(  \sum_{\ell=\mu}^{-\mu}\alpha_{\ell}^{0}y^{\ell} - \sum_{\ell=\mu}^{-\mu}\alpha_{-\ell}^{0}y^{\ell}
\right).
\ee
Comparing Eq. (\ref{odd2})  with Eq. (\ref{odd}) implies
\bee                     \nonumber
\alpha^0_{\ell}&=&\alpha^0_{-\ell}, \qquad \text{ if $\ell \ne \mu$, $\ell\ne -\mu$},  \\
\alpha^0_{\mu}&-&\alpha^0_{-\mu}= -\spur L_0,
\eee
and
\be\label{127}
F_{even}=\frac {\varkappa^{\mu}}{2}\alpha^0_{\mu}(y^{\mu}+y^{-\mu})
+\frac {\varkappa^{\mu}}{2} \sum_{\ell=0}^{\mu-1} \alpha^0_{\ell}(y^{\ell}+y^{-\ell})=
 \alpha^0_{\mu}\cosh(it\mu)
+ {\varkappa^{\mu}} \sum_{\ell=0}^{\mu-1} {\varkappa^{\ell}}\alpha^0_{\ell}\cosh (it\ell) .
\ee

\medskip

{\proposition
\be
\nonumber
\mathcal F^{\spur}(t) =
\alpha^0_{\mu}
+ {\varkappa^{\mu}} \sum_{\ell=0}^{\mu-1} {\varkappa^{\ell}}\alpha^0_{\ell}\cosh \left(t\sqrt{\mu^2-\ell^2}\,\right) .
\ee
}

\medskip

\medskip

\begin{proof}
Taking Proposition \ref {p121}  into account let us decompose Eq (\ref{121}) into the Taylor series:
\be            \nonumber
\mathcal F^{\spur}(t)=\sum_{s=0}^{\infty}a_{2s}\frac {t^{2s}}{(2s)!},
\ee
where $a_{2s}:=\spur(\m^{2s}Q_0)$ for $s=0,1,2,\,...$.

Eq (\ref{even}) implies
\be                 \nonumber
a_{2s}=\left( \frac {d^2}{dt^2} + \mu^2\right)^s F_{even}|_{t=0},
\ee
and Eq (\ref{127}) implies
\be                   \nonumber
a_{2s}=\left\{
\begin{array}{c}
a_{\mu }^{0} +
{\varkappa^{\mu}} \sum_{\ell=0}^{\mu-1} {\varkappa^{\ell}}\alpha^0_{\ell}
\text{ \ \ \ if }s=0 \\
{\varkappa^{\mu}} \sum_{\ell=0}^{\mu-1} {\varkappa^{\ell}}\alpha^0_{\ell}(-\ell^2+\mu^2)^s
 \text{ \ \ \  if }s\ne0    .
\end{array}%
\right.
\ee
So
\be                       \nonumber
\mathcal F^{\spur}(t)=\sum_{s=0}^{\infty}a_{2s}\frac {t^{2s}}{(2s)!}=
\alpha^0_{\mu}
+ {\varkappa^{\mu}} \sum_{\ell=0}^{\mu-1} {\varkappa^{\ell}}\alpha^0_{\ell}\cosh \left(t\sqrt{\mu^2-\ell^2}\,\right)  .
\ee
\end{proof}

\section{Ideals generated by degenerate $\varkappa$-traces}

{Let $\mu=z \in\mathbb Z \setminus n\mathbb Z$ and the $\varkappa$-trace be defined by
Eq (\ref{913}) for the case $\varkappa=+1$  and
by Eq (\ref{914}) for the case $\varkappa=-1$.

These degenerate $\varkappa$-traces are denoted in Theorem \ref{th2} by $tr_z$ and
$str_z$.

Denote the ideals generated by these traces $\spur$ by $\mathcal I^{\varkappa}$;   in $H^0$, consider
the ideals  $\mathcal I_0^{\varkappa} := \mathcal I^{\varkappa}\bigcap H^0$.

Now we can prove Conjecture \ref{con1} (\cite[Conjecture 9.1]{KT7}):

{\theorem\label{conj9.1}
$\mathcal I^{+1}=\mathcal I^{-1}$.
}

To prove Theorem \ref{conj9.1} we use Theorem 4.2 from \cite{KT20} which  in our case implies

{\theorem\label{4.2} (\cite[Theorem 4.2]{KT20})
$\mathcal I^{+1}=\mathcal I^{-1}$ if and only if $\mathcal I^{+1}\bigcap H^0=\mathcal I^{-1}\bigcap H^0$.
}

So, Theorem \ref{conj9.1} follows from 
{\theorem\label{4.3}
$\mathcal I^{+1}\bigcap H^0=\mathcal I^{-1}\bigcap H^0$.
}

\begin{proof}
For degenerate $\spur$, we established the following facts:
\bee             \nonumber
&&F_p^{\spur}(t)=\spur\left(e^{t\m}Q_p \right)=\varkappa^\mu\sum_{\ell=\mu}^{-\mu}
\alpha_{\ell}^{p}\varkappa^{\ell}
e^{it\ell} \text{ \ for \ } p=1,2,...,n-1  ,\\
&&\mathcal F^{\spur}(t)=\spur\left(e^{t\m}Q_0   \right)=
\alpha^0_{\mu}
+ {\varkappa^{\mu}} \sum_{\ell=0}^{\mu-1} {\varkappa^{\ell}}\alpha^0_{\ell}\cosh \left(t\sqrt{\mu^2-\ell^2}\,\right)
\text{\ \ where\ \ } \alpha^0_{\mu}\ne 0,
\nonumber
\eee
and where the $\alpha$-s do not depend on $\varkappa$.


For any $p=1,...,n$, it is easy to find the lowest degree polynomial differential operators with constant
coefficients $D^{\varkappa}_p(d/dt)$
 such that $D^{\varkappa}_p(d/dt) F_p^{\spur}(t)=0$:
\be              \nonumber
D^{\varkappa}_p\left(\frac d{dt}\right)= \begin{cases}
\mathop{\prod}\limits_{\ell=-\mu:\ \alpha_{\ell}^{p}\ne 0}^{\mu}
\left(\frac d {dt} -i\ell \right)&\text{\ \ if $F_p^{\spur}\ne 0$,}\\ 1&
\text{\ \ if $F_p^{\spur} = 0$,}\end{cases}
\ee
and $D^{\varkappa}_0(d/dt)$ such that $D^{\varkappa}_0(d/dt) \mathcal F^{\spur}(t)=0$:
\be                \nonumber
D^{\varkappa}_0\left(\frac d{dt}\right)=\frac d {dt}
\prod_{\ell=0:\ \alpha_{\ell}^{0}\ne 0}^{\mu-1}
\left(\frac {d^2 }{dt^2} -\mu^2+ \ell^2 \right).
\ee
Further, it is a~simple exercise to prove that
\be             \nonumber
D_p^{\varkappa}(\m)Q_p,\ D_p^{\varkappa}(\m)L_{-p} \in \mathcal I^{\varkappa}_0
\text{\ \ \ for any\ \ \ } p=0,\,...\,n-1,
\ee
namely,
\be               \nonumber
B_{\spur}(D_p^{\varkappa}(\m)Q_p,\, f )=B_{\spur}(D_p^{\varkappa}(\m)L_{-p},\, f )=0
\text{\ \ \ for any $f\in H^0$ and $p=0,\,...\,n-1$}.
\ee

Consider, for example, $B_{\spur}(D_0^{\varkappa}(\m)Q_0,\, f)$ for $f=g(\m)Q_p$ and $f=g(\m)L_p$:
\be             \nonumber
\spur(D_0^{\varkappa}(\m)Q_0 g(\m)Q_p)=\spur(D_0^{\varkappa}(\m)Q_0 g(\m)L_p)=0 \text{\ \ for \ \ }p\ne0,
\ee
since $Q_0 Q_p=Q_0 L_p=0$ for $p\ne 0$,
\bee
\spur(D_0^{\varkappa}(\m)Q_0 g(\m)Q_0)=\spur(D_0^{\varkappa}(\m) g(\m)Q_0)&=&
D_0^{\varkappa}\left(\frac d {dt}\right)g\left(\frac d {dt}\right) \spur\left(e^{t\m}Q_0    \right)|_{t=0}=  \nn
&=&g\left(\frac d {dt}\right) D_0^{\varkappa}\left(\frac d {dt}\right) \mathcal F^{\spur}(t)|_{t=0}=0,
\nonumber
\eee
\bee
\spur(D_0^{\varkappa}(\m)Q_0 g(\m)L_0) = \spur(D_0^{\varkappa}(\m) g(\m)L_0)&=&
D_0^{\varkappa}\left(\frac d {dt}\right)g\left(\frac d {dt}\right) \spur\left(e^{t\m} L_0   \right)|_{t=0}=  \nn
&=&g\left(\frac d {dt}\right) D_0^{\varkappa}\left(\frac d {dt}\right) \spur(L_0)=0
\nonumber
\eee
due to Eq (\ref{psiq}) and since the operator $D_0^{\varkappa}\left(\frac d {dt}\right)$
contains the factor $\frac d {dt}$.


Further, it is easy to see that for each of the ideals $\mathcal I_0^{\varkappa}$, where $\varkappa=\pm 1$,
the polynomials $\phi_p^0\in \mathbb C[\m]$ defined in Corollary \ref{col} satisfy the relations
$\phi_p^0(\m)=D_p^{\varkappa}(\m)$ for $p=0,...,n-1$.

So, Theorem \ref{i0} implies that the $\mathbb C [\m]$-span of the $D_p^{\varkappa}(\m)Q_p$ and $D_p^{\varkappa}(\m)L_{-p}$
for $p=0,...,n-1$
is $\mathcal I^{\varkappa}_0$.

Since $D^{+1}_p=D^{-1}_p$, we have $\mathcal I^{+1}_0=\mathcal I^{-1}_0$,
and as result, $\mathcal I^{+1}=\mathcal I^{-1}$.
\end{proof}

\vskip 5mm

\section*{Acknowledgments}
The authors
are grateful to Russian Fund for Basic Research
(grant No.~${\text{20-02-00193}}$)
for partial support of this work.


\setcounter{equation}{0} \def\theequation{A
\arabic{equation}}

\newcounter{appen}
\newcommand{\appen}[1]{\par\refstepcounter{appen}
{\par\bigskip\noindent\large\bf Appendix. 
\medskip }{\bf \large{#1}}}

\renewcommand{\subsection}[1]{\refstepcounter{subsection}
\vskip 3mm{\centerline{\bf 
A\arabic{subsection}. \ #1}}\vskip 3mm}
\renewcommand\thesubsection{A\theappen.\arabic{subsection}}
\makeatletter \@addtoreset{subsection}{appen}

\renewcommand{\subsubsection}{\par\refstepcounter{subsubsection}
{\bf A\arabic{appen}.\arabic{subsection}.\arabic{subsubsection}. }}
\renewcommand\thesubsubsection{A\theappen.\arabic{subsection}.\arabic{subsubsection}}
\makeatletter \@addtoreset{subsubsection}{subsection}


\newcounter{theor}
\renewcommand{\theorem}{\par\refstepcounter{theor}
           {\bf Theorem A.%
           \arabic{theor}. }}
\renewcommand\thetheorem{A.\arabic{theor}}
\makeatletter \@addtoreset{theor}{section}

\renewcommand{\remark}{\par\refstepcounter{theor}
           {\bf Remark A.%
           \arabic{theor}. }}
\renewcommand\theremark{\thesection.%
     \arabic{theor}}
\makeatletter \@addtoreset{remark}{section}

\renewcommand{\definition}{\par\refstepcounter{theor}
           {\bf Definition A.%
           \arabic{theor}. }}
\renewcommand\thedefinition{\thesection.%
    \arabic{theor}}
\makeatletter \@addtoreset{definition}{section}

\renewcommand{\proposition}{\par\refstepcounter{theor}
           {\bf Proposition 
           A.\arabic{theor}. }}
\renewcommand\theproposition{
  A.\arabic{theor}}
\makeatletter \@addtoreset{proposition}{section}

\renewcommand{\theorem}{\par\refstepcounter{theor}
           {\bf Theorem 
           A.\arabic{theor}. }}
\renewcommand\thetheorem{
   A.\arabic{theor}}
\makeatletter \@addtoreset{theorem}{section}




\end{document}